\documentclass[11pt]{article}
\usepackage[utf8]{inputenc}
\IfFileExists{lmodern.sty}{\usepackage[T1]{fontenc}\usepackage{lmodern}}{\usepackage[OT1]{fontenc}}
\usepackage{microtype}
\usepackage[letterpaper,hmargin=1in,top=0.83in,bottom=0.87in]{geometry}
\usepackage{amsmath,amssymb,amsthm,enumitem}
\usepackage[dvipsnames]{xcolor}
\usepackage[colorlinks=true,linkcolor=MidnightBlue,citecolor=MidnightBlue,urlcolor=MidnightBlue]{hyperref}
\hypersetup{pdftitle={Prescribed matrix size for UHF Brauer classes on suspensions},pdfauthor={Andrei V. Ershov},pdfsubject={Matrix embeddings, localized K-theory and unital representatives}}
\setlist[enumerate]{leftmargin=*,topsep=4pt,itemsep=2pt,parsep=0pt}
\newtheorem{theorem}{Theorem}[section]
\newtheorem{corollary}[theorem]{Corollary}
\newtheorem{proposition}[theorem]{Proposition}
\newtheorem{lemma}[theorem]{Lemma}
\DeclareMathOperator{\Aut}{Aut}
\DeclareMathOperator{\End}{End}
\DeclareMathOperator{\Emb}{Emb}
\DeclareMathOperator{\Fr}{Fr}
\DeclareMathOperator{\Br}{Br}
\newcommand{\K}{\mathcal K}
\newcommand{\Z}{\mathbb Z}
\newcommand{\C}{\mathbb C}
\newcommand{\id}{\mathrm{id}}
\newcommand{\triv}{\mathbf 1}
\newcommand{\RP}{\mathbb{RP}}
\newcommand{\RY}{R_Y}
\newcommand{\IY}{I_Y}
\title{Prescribed matrix size for UHF Brauer classes\\ on suspensions}
\author{Andrei V. Ershov}
\date{22 November 2026}
\begin{document}
\maketitle
\begin{abstract}
Let \(Y\) be a connected finite CW complex, \(L=M_{l^\infty}\) for \(l\geq2\), and \(k\) a positive integer coprime to \(l\). Stable \(L\)-twists on \(\Sigma Y\) have normalized clutching invariants \(1+x\) in \(K^0(Y)[1/l]\). We prove that such a twist is the stabilization of a unital \(M_k(L)\)-bundle exactly when \(kx=0\) additively. For each admissible size we construct a representative from two trivial \(M_k(\C)\)-bundles and a span of embeddings into \(M_{kl^N}(\C)\) on their overlap. Every \(N\geq0\) with \(2l^N\geq\dim Y\) suffices. For \(Y=S^{2m-1}\cup_p e^{2m}\), \(m\geq2\), and a prime \(p\geq m\) coprime to \(l\), a generator has least unital size \(p\) and least overlap exponent \(\lceil\log_l m\rceil\) at that size. In dimension five the least size is \(q\) and the least exponent is one for every Moore space \(S^4\cup_q e^5\) with \(q>1\) coprime to \(l\). Over \(\Sigma\RP^{2d}\), a twist of period two has least unital size \(2^d\) whenever \(l\) is odd.
\end{abstract}

\section{The fixed-size question}\label{sec:intro}

Let \(L=M_{l^\infty}=\bigotimes_{n\geq1}M_l(\C)\), let \(\K\) denote the compact operators on a separable infinite-dimensional Hilbert space, and write \(M_k(L)=M_k(\C)\otimes L\). A unital \(M_k(L)\)-bundle \(E\) yields a stable twist by taking \(E\otimes\K\). For connected compact metrizable \(X\), Dadarlat and Pennig identify \(\Br_L(X)\) with the torsion subgroup of the orientable stable twists \cite[Theorem~2.15(ii)]{DP2}. On the suspensions considered here the equator \(Y\) is connected, so \(X=\Sigma Y\) is simply connected. The orientation obstruction \(\delta_0\in H^1(X;\Z[1/l]^\times_+)\) therefore vanishes for every stable twist \cite[Proposition~2.3]{DP2}. For \(k\geq2\), Dadarlat and Pennig also show that an orientable stable twist annihilated by a power of \(k\) has a unital representative of size \(k^r\) for some \(r\) \cite[Theorem~2.11 and Corollary~2.12]{DP2}. Which prescribed sizes work for a given twist?

On a suspension \(X=\Sigma Y\), the invariant that answers this question is an element of the ordinary complex \(K\)-ring of the equator. Fix a basepoint in \(Y\) and set
\[
  \RY=K^0(Y)\otimes\Z[1/l],\qquad
  \IY=\widetilde K^0(Y)\otimes\Z[1/l].
\]
Throughout, \(Y\) is a connected finite CW complex. We use the unreduced suspension with its two cones. The reduced suspension of a pointed CW model is homotopy equivalent to it. A stable twist means a locally trivial \(L\otimes\K\)-bundle up to bundle isomorphism. We use the suspension clutching construction of \cite[Section~4.2]{Pennig}. The evaluation theorem \cite[Theorem~2.22]{DP1}, followed by the quotient that forgets the cone trivializations, identifies a stable twist with a normalized unit \(1+x\in 1+\IY\) (Section~\ref{sec:clutching}). The finite data for a two-patch lax algebra bundle (LAB) are two embeddings from the cone algebras into a common matrix algebra on the equator.

\newpage
\begin{theorem}[Prescribed size and finite overlap]\label{thm:main}
Let \(l\geq2\), \(k\geq1\), and \(\gcd(k,l)=1\). For a stable twist on \(\Sigma Y\) with normalized invariant \(1+x\), the following conditions are equivalent.
\begin{enumerate}[label=\textup{(\roman*)}]
\item The twist is isomorphic to \(E\otimes\K\) for a locally trivial unital \(M_k(L)\)-bundle \(E\).
\item \(kx=0\) in the additive group \(\IY\).
\item The twist is represented by a two-patch LAB of size \(k\) with one standard embedding on the overlap, at some finite stage.
\end{enumerate}
When these conditions hold, one may prescribe \emph{any} stage \(N\geq0\) such that \(2l^N\geq\dim Y\).
\end{theorem}

Thus the least admissible matrix size is the additive order of \(x\) if that order is finite, with size one for \(x=0\). The period of the twist is the multiplicative order of \(1+x\). These two numbers can differ. Prescribing \(k\) asks for the image of a stabilization map on bundle classes. Different unital representatives can have the same stable class.

The construction uses the multiplicity bundle of a family of matrix embeddings. Its role in the action of the frame space on \(K\)-theory appears in \cite[Proposition 30 and Theorem 33]{ErshovEmb}. Spans of such embeddings occur in the theory of lax algebra bundles \cite[Example~3.6]{ErshovLAB}. Here the two cones suffice to pass from that finite geometry to the stable invariant, independently of a classification theorem for general covers.

\section{Stable clutching and the unit corner}\label{sec:clutching}

All automorphism groups carry the point-norm topology. Set \(G_{\mathrm{st}}=\Aut(L\otimes\K)\) and choose a rank-one projection \(e\in\K\). The evaluation isomorphism of \cite[Theorem 2.22]{DP1} sends a continuous family \(\gamma:Y\to G_{\mathrm{st}}\) to
\begin{equation}\label{eq:evaluation}
    [\,y\longmapsto\gamma_y(1_L\otimes e)\,]
       \in K_0(C(Y)\otimes L)^\times_+ .
\end{equation}
Composition of families induces multiplication of units. The identification \(K_0(C(Y)\otimes L)=\RY\) is normalized by \([1_L]=1\). Indeed, under matrix Morita equivalence, the transition on \(K^0(Y)\) from the \(M_{l^N}\) stage to the \(M_{l^{N+1}}\) stage is multiplication by \(l\). A projection with range bundle \(\eta\) at stage \(N\) therefore contributes \([\eta]/l^N\).

To clutch a bundle, choose trivializations on the cones and take the transition from the negative cone to the positive cone. Changing either trivialization multiplies its class in \([Y,G_{\mathrm{st}}]\) by a constant family, since a map on a cone is homotopic to a constant. Conversely, these changes and homotopies of the transition give isomorphic bundles. The group \([Y,G_{\mathrm{st}}]\) is abelian by \eqref{eq:evaluation}. Thus the quotient by the constant families is computed by dividing the evaluation unit by its rank at the basepoint. The constant positive units are exactly \(\Z[1/l]^\times_+\).

To see that the normalized positive units are all of \(1+\IY\), take \(v=1+x\) with \(x\in\IY\). The skeletal filtration makes \(\IY\) nilpotent, so \(v^{-1}\) exists and also has rank one. For either \(w=v\) or \(w=v^{-1}\), choose \(m\) large enough that \(l^mw\) has an integral lift in \(K^0(Y)\) and \(2l^m\geq\dim Y\). The lift has rank \(l^m\), so the complex stable-range theorem used in the proof of Theorem~\ref{thm:main} represents it by an actual vector bundle. Represent that bundle by a projection in a matrix amplification of the stage algebra. In the stage-\(m\) normalization the projection has class \(w\). This gives positivity for both \(v\) and its inverse. In reduced coordinates the tensor law is \(x\star y=x+y+xy\), since \((1+x)(1+y)=1+x+y+xy\).

Write \(D_k=M_k(L)\). To fix \(D_k\otimes\K\cong L\otimes\K\), absorb the \(M_k\) factor into \(\K\) using a Hilbert-space unitary that takes the rank-one projection \(e_{11}\otimes e\) to \(e\).

\begin{lemma}\label{lem:corner}
Let \(\alpha:Y\to\Aut(D_k)\) clutch a unital bundle. In the Morita identification \(K_0(C(Y)\otimes D_k)=\RY\), use the normalization \([e_{11}\otimes1_L]=1\). The stabilized bundle has normalized invariant
\[
  u=[\,y\longmapsto\alpha_y(e_{11}\otimes1_L)\,],
  \qquad ku=k .
\]
\end{lemma}
\begin{proof}
Formula~\eqref{eq:evaluation} applied to \(\alpha\otimes\id_{\K}\) gives the stated projection class. The \(k\) families \(\alpha(e_{ii}\otimes1_L)\) are mutually equivalent via the partial isometries \(\alpha(e_{ij}\otimes1_L)\). Their sum is the constant unit of \(D_k\), of class \(k\). Hence \(ku=k\). At a point this equation holds in the torsion-free group \(\Z[1/l]\), so the rank of \(u\) is already one.
\end{proof}

\section{Realization by a finite span}\label{sec:frames}

Write \(M_t=M_t(\C)\). A two-patch LAB of size \(k\) and stage \(N\) consists here of a trivial \(M_k\)-bundle on each cone, a trivial \(M_{kl^N}\)-bundle on the equator, and a continuous unital embedding from each restricted cone algebra into the equatorial algebra. Pullback over an equatorial collar gives the corresponding span on open cone neighborhoods. For a stable twist and a fixed cone size \(k\), its \emph{least overlap exponent at size \(k\)} is the smallest \(N\geq0\) for which such a span represents the twist after stabilization, when such a span exists. The minimum ranges over both embedding legs, including spans in which neither leg is standard.

For \(t=l^N\), let \(\Fr_{k,t}=\Emb(M_k,M_{kt})\) be the space of unital \(*\)-embeddings. Its topology is determined by norm convergence on the matrix units. The standard embedding is \(\iota(a)=a\otimes 1_t\).

\begin{lemma}[Multiplicity bundle]\label{lem:mult}
A continuous map \(f:Y\to\Fr_{k,t}\) determines a rank-\(t\) bundle \(\eta_f\) and a unitary trivialization of \(\C^k\otimes\eta_f\). Conversely, a rank-\(t\) Hermitian bundle \(\eta\) equipped with such a trivialization gives a frame family \(f\) with \(\eta_f\cong\eta\). The commutant bundle of \(f(M_k)\) in the trivial \(M_{kt}\)-bundle is canonically \(\End(\eta_f)\).
\end{lemma}
\begin{proof}
Take \((\eta_f)_y=f_y(e_{11})\C^{kt}\). These projections vary continuously and have rank \(t\). If \(v_1,\ldots,v_k\) is the standard basis of \(\C^k\), the rule \(v_i\otimes w\mapsto f_y(e_{i1})w\) gives the trivialization. It identifies \(f_y(M_k)\) with the operators on the first tensor factor, so the commutant acts on the second. Conversely, transport the action \(a\mapsto a\otimes1_\eta\) through the given trivialization.
\end{proof}

\begin{proposition}[Clutching formula for a span]\label{prop:span}
For a span \(M_k\xrightarrow{f_-}M_{kt}\xleftarrow{f_+}M_k\) parametrized by \(Y\), there is a well-defined isomorphism class of unital \(D_k\)-bundles on \(\Sigma Y\). Its stabilized clutching invariant is
\begin{equation}\label{eq:span}
  u=\bigl([\eta_{f_+}]/t\bigr)^{-1}\bigl([\eta_{f_-}]/t\bigr).
\end{equation}
In particular, for \(f_+=\iota\) the invariant is \([\eta_{f_-}]/t\).
\end{proposition}
\begin{proof}
For either frame \(f\), the commutant after tensoring the ambient algebra by \(L\) is \(\End(\eta_f)\otimes L\). It is a locally trivial bundle with fiber \(M_t\otimes L\cong L\). Since \(Y\) is paracompact and \(\Aut(L)\) is contractible \cite[Theorem 2.3]{DP1}, it admits a trivialization and any two trivializations are homotopic. Choose one in the direction \(L\to\End(\eta_f)_y\otimes L\). Multiplication with \(f_y(M_k)\otimes1\) now gives a continuous family of isomorphisms
\[
 \widetilde f_y:D_k\longrightarrow M_{kt}\otimes L,\qquad
 \widetilde f_y(a\otimes1)=f_y(a)\otimes1 .
\]
Fix the tail identification \(\theta:M_{kt}\otimes L\to D_k\) induced by \(M_{kt}=M_k\otimes M_t\) and \(M_t\otimes L\cong L\) for \(t=l^N\). Set \(\alpha_f=\theta\circ\widetilde f:Y\to\Aut(D_k)\), where composition is fiberwise. For the standard frame take \(\widetilde\iota=\theta^{-1}\), so \(\alpha_\iota=\id\). The transition from the negative cone to the positive cone is
\[
 \gamma=\alpha_{f_+}^{-1}\alpha_{f_-}
        =\widetilde f_+^{-1}\widetilde f_- .
\]
The choices of commutant trivializations vary through homotopies, so \(\gamma\) determines an isomorphism class of unital bundles.

Use matrix Morita coordinates on the source and target of \(\widetilde f\), in each case assigning class \(1\) to a rank-one matrix corner tensored with \(1_L\). With \(p=e_{11}\otimes1_L\) in \(D_k\), the image \(\widetilde f(p)=f(e_{11})\otimes1_L\) has class \([\eta_f]\) in the target. The induced map on \(K_0\) commutes with tensoring by bundles on \(Y\) because the family is \(C(Y)\)-linear. It is therefore \(K^0(Y)\)-linear, and hence \(R_Y\)-linear after inverting \(l\). Its value on \(1=[p]\) shows that it acts by multiplication by \([\eta_f]\) in these coordinates.

Under \(\theta\), the constant rank-\(t\) corner \(e_{11}\otimes1_t\otimes1_L\) goes to \(p\). Its \(t\) rank-one summands are equivalent, so \(\theta_*\) sends each summand to \(1/t\). Applying Lemma~\ref{lem:corner} now gives \(u(\alpha_f)=[\eta_f]/t\). Multiplicativity of \eqref{eq:evaluation} gives \(u(\gamma)=u(\alpha_{f_+})^{-1}u(\alpha_{f_-})\), as claimed.
\end{proof}

\begin{lemma}\label{lem:torsion}
For any abelian group \(A\) and \((k,l)=1\), localization induces an isomorphism \(A[k]\to A[1/l][k]\), where \(A[k]=\ker(k:A\to A)\).
\end{lemma}
\begin{proof}
Flatness of \(\Z[1/l]\) identifies the kernel after localization with \(A[k]\otimes\Z[1/l]\). Multiplication by \(l\) is already invertible on \(A[k]\), by B\'ezout's identity.
\end{proof}

\begin{proof}[Proof of Theorem~\ref{thm:main}]
Lemma~\ref{lem:corner} gives \textup{(i)}\(\Rightarrow\)\textup{(ii)}, since \(k(1+x)=k\). Proposition~\ref{prop:span} gives \textup{(iii)}\(\Rightarrow\)\textup{(i)}.

Suppose \(kx=0\). By Lemma~\ref{lem:torsion}, \(x\) has a unique lift \(a\in\widetilde K^0(Y)[k]\). Fix \(t=l^N\) with \(2t\geq\dim Y\). In integral \(K^0(Y)\), form the rank-\(t\) class \(v=t[\triv]+t a\), where \(\triv\) is the trivial complex line. Complex vector bundles of rank \(s\) are classified by their stable class on CW complexes of dimension at most \(2s\). Indeed \(BU(s)\to BU\) is \((2s+1)\)-connected, as follows from \(U(s+1)/U(s)=S^{2s+1}\) and stabilization. Thus \(v\) is represented by a rank-\(t\) bundle \(\eta\). As \(ka=0\), the \(k\)-fold direct sum of \(\eta\) has stable class \(kt[\triv]\). The same range applied at rank \(kt\) gives an actual isomorphism
\[
                    \eta^{\oplus k}\ \cong\ Y\times\C^{kt}.
\]
Choose a unitary trivialization. Lemma~\ref{lem:mult} supplies \(f:Y\to\Fr_{k,t}\) with multiplicity bundle \(\eta\). Take \(f_- = f\) and \(f_+=\iota\). Formula~\eqref{eq:span} gives \(u=[\eta]/t=1+x\). This proves \textup{(iii)} with the asserted stage bound.
\end{proof}

\begin{corollary}\label{cor:orders}
If \(x\) has finite additive order \(a\), a positive integer \(s\) occurs as the unital size of a representative of \(1+x\) exactly when \(a\mid s\). If \(x\) has infinite additive order there is no unital representative of any size. As \(k\) varies, the stable classes obtained are exactly the torsion elements of \(1+\IY\).
\end{corollary}
\begin{proof}
Lemma~\ref{lem:corner} gives \(sx=0\) for every unital size \(s\), without any coprimality assumption. A finite additive order \(a\) in \(\IY\) is prime to \(l\), since multiplication by each prime dividing \(l\) is invertible. Theorem~\ref{thm:main} therefore supplies a representative of size \(a\). Tensoring it with a constant matrix algebra realizes each multiple of \(a\) with the same stabilized bundle, since \(M_b\otimes\K\cong\K\) for every \(b\geq1\). If \(x\) has finite additive order, the powers \(x,x^2,\ldots\) generate a finite additive group, since each is torsion of bounded order and \(\IY\) is nilpotent. The powers of \(1+x\) lie in a finite set and \(1+x\) is a unit. Conversely, if \((1+x)^m=1\) for some \(m\geq1\), then
\[
 x\left(\sum_{j=1}^{m}\binom{m}{j}x^{j-1}\right)=0.
\]
The sum has constant term \(m\) and nilpotent remainder, hence is a unit over \(\mathbb Q\). Thus \(x\) is additively torsion.
\end{proof}

\section{Two examples}\label{sec:examples}

\subsection{A higher Brauer class requiring one enlargement}

\begin{theorem}\label{thm:moore}
Let \(q>1\), \((q,l)=1\), \(Y=S^3\cup_q e^4\), and \(X=\Sigma Y\simeq S^4\cup_q e^5\), where the attaching maps have degree \(q\). Then \(\Br_L(X)\cong\Z/q\). A generator has least unital size \(q\) and is realized by a span
\[
                M_q\longrightarrow M_{ql}\longleftarrow M_q .
\]
Its least overlap exponent in the \(l\)-power tower at that size is \(N=1\). No ordinary finite matrix algebra bundle of any size yields this generator by extension to \(L\) and stabilization.
\end{theorem}
\begin{proof}
The reduced \(K\)-theory exact sequence of \(S^3\xrightarrow{q}S^3\to Y\) identifies \(\widetilde K^0(Y)\) with the cokernel of \(q:\widetilde K^{-1}(S^3)\to\widetilde K^{-1}(S^3)\), hence with \(\Z/q\) by Bott periodicity. The space \(Y\) is a suspension up to homotopy, so its reduced diagonal is null-homotopic. Products of reduced \(K^0\)-classes therefore vanish. Localization at \(l\) preserves \(\Z/q\). It follows that \(1+\IY\) is cyclic of order \(q\) under multiplication. Every such class belongs to \(\Br_L(X)\) by \cite[Theorem~2.15(ii)]{DP2}. A generator \(1+x\) has additive coordinate of order \(q\), so its least size is \(q\).

The bound \(2l\geq4=\dim Y\) in Theorem~\ref{thm:main} gives the displayed span at \(N=1\). More explicitly, take a rank-\(l\) bundle of class \(l[\triv]+l a\), where \(a\) is the integral lift of \(x\), and trivialize its \(q\)-fold direct sum. This produces the varying leg of the span.

At \(N=0\) the two legs are isomorphisms of \(M_q\). Composing one with the inverse of the other gives an ordinary \(M_q\)-clutching map. But \(H^3(X;\Z)=0\) by the cellular cochain complex. For any ordinary \(M_s\)-bundle \(B\) on \(X\), the obstruction to lifting its principal \(\mathrm{PU}(s)\)-bundle to \(U(s)\) lies in this group and hence vanishes. Thus \(B\cong\End(V)\) for a vector bundle \(V\), and \(B\otimes\K\) has zero Dixmier--Douady class. Consequently \(B\otimes L\otimes\K\) is trivial. Since \(1+x\) is nontrivial, every \(N=0\) presentation is excluded, even if the cone size is allowed to vary.
\end{proof}

Here a higher class is detected by its failure to lie in the image of the ordinary topological Brauer group. Dadarlat and Pennig use the same five-dimensional Moore space for \(q=3\), with Jiang--Su coefficients, to study the opposite-algebra operation \cite[Example 3.5]{DP2}. Theorem~\ref{thm:moore} determines the sizes of the two finite matrix algebras that exhibit a generator.

The varying leg also detects the tautological example from the frame space \cite[Example~3.6]{ErshovLAB}. On \(\Sigma\Fr_{q,l}\) use the identity family of frames on one leg and the standard embedding on the other. For the map \(f:Y\to\Fr_{q,l}\) constructed above, pullback of this span along \(\Sigma f\) gives the Moore-space generator. If the stabilized tautological class came from any ordinary finite matrix bundle, its pullback would have the same property. Thus the finite frame space already carries, after suspension, a Brauer class outside the classical image.

\paragraph{The twisted \(K\)-theory rank map.}
For a generator \(\tau\) in Theorem~\ref{thm:moore}, write \(1+g\) for
its clutching unit and \(R=\Z[1/l]\). Write \(KU[1/l]^i(X;\tau)\)
for the \(K_i\)-group of the section algebra of the stable bundle
representing \(\tau\). Pennig's suspension calculation
\cite[Theorem~4.5]{Pennig} identifies evaluation at a cone point
in the Mayer--Vietoris sequence. Here \(K^1(Y)[1/l]=0\) by the
cofiber sequence, and \(K^0(Y)[1/l]=R\oplus(R/qR)g\). With the two
cone trivializations, the relevant exact sequence is
\[
 0\longrightarrow KU[1/l]^0(X;\tau)\longrightarrow R\oplus R
 \xrightarrow{\ \delta\ } R\oplus(R/qR)g
 \longrightarrow KU[1/l]^1(X;\tau)\longrightarrow 0 .
\]
In coordinates ordered so that the second cone is transported to
the first, \(\delta(a,b)=a-(1+g)b=(a-b,-bg)\).
The map is surjective, and its kernel consists of the pairs
\((a,a)\) with \(a\in qR\). Projection to either cone coordinate
therefore gives
\[
 KU[1/l]^0(X;\tau)\xrightarrow[\cong]{\ \mathrm{ev}\ }qR\subset R,
 \qquad KU[1/l]^1(X;\tau)=0 .
\]
Evaluation records the divisibility of fiber ranks by \(q\).
The unit of the size-\(q\) representative in Theorem~\ref{thm:moore}
attains the positive integer rank \(q\).

\subsection{Period two and unbounded least size}

\begin{proposition}\label{prop:projective}
Let \(l\) be odd, \(d\geq1\), and \(Y=\RP^{2d}\). Let \(\lambda\) be the complexification of the tautological real line. The stable twist on \(\Sigma Y\) with clutching unit \([\lambda]\) has period \(2\) and least unital size \(2^d\). At that size it already has a presentation at \(N=0\).
\end{proposition}
\begin{proof}
The standard complex \(K\)-theory calculation gives
\(\widetilde K^0(\RP^{2d})=\Z/2^d\), generated by \(a=[\lambda]-1\)
\cite[Theorem 4.6.1]{RW}. Its order survives localization at any odd \(l\). Since \(\lambda^{\otimes2}\cong\triv\), the unit \([\lambda]\) has order two. Its least size follows from Corollary~\ref{cor:orders}.

Put \(k=2^d\). Integrally \(k[\lambda]=k[\triv]\), and \(2d\leq2k\). Cancellation in rank \(k\) gives a trivialization of \(\lambda^{\oplus k}\). Apply Lemma~\ref{lem:mult} with \(t=1\) and \(\eta=\lambda\). It gives a map \(Y\to\Aut(M_k)\). The span with its other leg the identity therefore clutches an ordinary \(M_k\)-bundle with the required invariant.
\end{proof}

In this example the minimal representative is already an ordinary matrix bundle. It illustrates how a classical period--index gap persists after adjoining an odd UHF tail. For the general topological period--index problem for ordinary matrix bundles, see Antieau and Williams \cite{AW}. In Theorem~\ref{thm:moore}, the period and least size agree, yet a finite enlargement on the overlap is necessary. Thus the least unital size and the least overlap stage measure different aspects of the finite presentation.

\section{Sharp overlap exponents on higher Moore spaces}\label{sec:sharp-overlap}

The dimension bound in Theorem~\ref{thm:main} is attained in
arbitrarily high dimensions. The obstruction already appears in the
top Chern class of the multiplicity bundle of each embedding leg.

\begin{theorem}\label{thm:sharp-moore}
Let \(m\geq2\), let \(p\geq m\) be prime, and assume \(p\nmid l\).
Put \(Y=S^{2m-1}\cup_p e^{2m}\) and \(X=\Sigma Y\). Then
\(\Br_L(X)\cong\Z/p\). Each generator has least unital size \(p\).
At cone size \(p\), its least overlap exponent is
\[
 N_{\min}=\min\{N\geq0:l^N\geq m\}=\lceil\log_l m\rceil .
\]
The minimum is taken over spans in which both embeddings may vary.
\end{theorem}
\begin{proof}
Let \(c:Y\to S^{2m}\) collapse the bottom sphere, and choose a Bott
generator \(b\in\widetilde K^0(S^{2m})\). The cofiber sequence of
the attaching map identifies \(c^*\) in reduced \(K^0\) with the
quotient \(\Z\to\Z/p\). Thus \(c^*b\) generates
\(\widetilde K^0(Y)\), and localization at \(l\) preserves its order.
The space \(Y\) is a suspension up to homotopy, so products in its
reduced \(K\)-ring vanish. Hence addition in \(\IY\) corresponds to
multiplication in \(1+\IY\). The clutching description and
\cite[Theorem~2.15(ii)]{DP2} give \(\Br_L(X)\cong\Z/p\).
Corollary~\ref{cor:orders} gives least size \(p\) for a generator.
Since \(\dim Y=2m\), Theorem~\ref{thm:main} realizes it at size \(p\)
whenever \(l^N\geq m\).

For the converse, all \(H^{2i}(Y;\Z)\) with \(0<i<m\) vanish.
The Whitney sum formula therefore makes the Chern class
\[
 c_m:\widetilde K^0(Y)\longrightarrow H^{2m}(Y;\Z)
\]
additive on virtual classes. Let \(z\) generate \(H^{2m}(S^{2m};\Z)\).
Bott integrality gives \(c_m(b)=\pm(m-1)!\,z\)
\cite[p.~308, Corollary~9.8]{Husemoller}. The cellular cochain
sequence identifies \(c^*z\) with a generator of \(H^{2m}(Y;\Z)=\Z/p\).
By naturality, \(c_m\) sends \(c^*b\) to \(\pm(m-1)!\,c^*z\).
Since \(p\geq m\), this is a generator, and \(c_m\) is injective.

Suppose now that \(t=l^N<m\). For any rank-\(t\) bundle \(\eta\) on
\(Y\), one has \(c_m([\eta]-t[\triv])=c_m(\eta)=0\), since \(m>t\).
Injectivity gives \([\eta]=t[\triv]\) in integral \(K^0(Y)\).
In particular the two multiplicity bundles of any span at stage \(N\)
both have class \(t[\triv]\). Formula~\eqref{eq:span} then gives
clutching unit \(1\). A generator consequently requires \(l^N\geq m\).
\end{proof}

For \(l=2\) and \(p=m=3\), the base is \(X\simeq S^6\cup_3 e^7\).
At the least cone size \(3\), the first overlap algebra that realizes
a generator is \(M_{12}\) at \(N=2\). An \(M_6\) overlap at \(N=1\)
has rank-two multiplicity bundles, whose third Chern classes vanish.
More generally, fixing \(l\) and taking \(m\) arbitrarily large with
a prime \(p\geq m\) coprime to \(l\) gives unbounded least overlap
exponents.

For a fixed prime \(p\nmid l\) and \(m>p\), the same calculation
gives the zero homomorphism \(c_m\), since \(p\mid(m-1)!\).
The upper bound of Theorem~\ref{thm:main} still applies. Determining
the least stage in this range requires information about which
reduced \(K\)-classes arise from bundles of small rank on \(Y\).

\clearpage
\noindent\textit{Use of generative AI.} OpenAI's ChatGPT assisted in developing and checking the arguments, searching the literature, and drafting the exposition and \LaTeX{} source. Responsibility for the mathematical verification, references, and final text rests with the author.

\end{document}